# Identities inspired from the Ramanujan Notebooks

# First series (1998)

by Simon Plouffe

## Abstract

I present here a collection of formulas inspired from the Ramanujan Notebooks. These formulas were found using an experimental method based on three widely available symbolic computation programs: PARI-Gp, Maple and Mathematica. A new formula is presented for $\zeta(5)$.

Une collection de formules inspirées des Notebooks de S. Ramanujan, elles ont toutes été trouvées par des méthodes expérimentales. Ces programmes de calcul symbolique sont largement disponibles (Pari-GP, Maple, Mathematica). Une nouvelle formule pour $\zeta(5)$ est présentée.

# Introduction

The first idea was to find the natural generalization of the $\zeta(3)$ result from Ramanujan, there are many possible generalizations to choose from, that formula appears in Ramanujan Notebooks II, chapter 14 formulas 25.1 and 25.3. Formula 25.1 [6].

$$\zeta(3) = \frac{7\pi^3}{180} - 2\sum_{n=1}^{\infty} \frac{1}{n^3(e^{2\pi n} - 1)}$$

A first step was to use a variation in the sign. The goal was to get an interesting identity for $\zeta(5)$.

$$\zeta(3) = 14\sum_{n=1}^{\infty} \frac{1}{n^3 \sinh(\pi n)} - \frac{11}{2}\sum_{n=1}^{\infty} \frac{1}{n^3(e^{2\pi n} - 1)} - \frac{7}{2}\sum_{n=1}^{\infty} \frac{1}{n^3(e^{2\pi n} + 1)}$$

$$\zeta(5) = \frac{\pi^5}{294} - \frac{72}{35}\sum_{n=1}^{\infty} \frac{1}{n^5(e^{2\pi n} - 1)} - \frac{2}{35}\sum_{n=1}^{\infty} \frac{1}{n^5(e^{2\pi n} + 1)}$$

These 2 first formulas found, it was not difficult to guess the next ones.

$$\zeta(5) = -\frac{39}{20}\sum_{n=1}^{\infty} \frac{1}{n^5(e^{2\pi n} - 1)} + \frac{1}{20}\sum_{n=1}^{\infty} \frac{1}{n^5(e^{2\pi n} + 1)} + 12\sum_{n=1}^{\infty} \frac{1}{n^5 \sinh(\pi n)}$$

Then for $\zeta(7)$

$$\zeta(7) = \frac{19\pi^7}{57600} - 2\sum_{n=1}^{\infty} \frac{1}{n^7(e^{2\pi n} - 1)}$$

The new formula being the first $\zeta(5)$ formula.

Another direction was suggested by the Apéry formula for $\zeta(3)$.

$$\zeta(3) = \frac{5}{2}\sum_{n=1}^{\infty} \frac{(-1)^{n+1}}{n^3 \binom{2n}{n}}$$

From which the case for $\zeta(5)$ doesn't seems to be there. One possible direction is to use the $\psi$ function and derivatives.

$$\sum_{n=1}^{\infty} \frac{1}{n^3 \binom{2n}{n}} = \frac{\pi\sqrt{3}}{18}\left(\psi'\left(\frac{1}{3}\right) - \psi'\left(\frac{2}{3}\right)\right) - \frac{4}{3}\zeta(3)$$

Which of course leads to $\zeta(5)$ and $\zeta(7)$.

$$\sum_{n=1}^{\infty} \frac{1}{n^5 \binom{2n}{n}} = \frac{\pi\sqrt{3}}{432}\left(\psi'''\left(\frac{1}{3}\right) - \psi'''\left(\frac{2}{3}\right)\right) + \frac{19}{5}\zeta(5) + \frac{\zeta(3)\pi^2}{9}$$

$$\sum_{n=1}^{\infty} \frac{1}{n^7 \binom{2n}{n}} = \frac{11\pi\sqrt{3}}{311040}\left(\psi'''''\left(\frac{1}{3}\right) - \psi'''''\left(\frac{2}{3}\right)\right) - \frac{493}{24}\zeta(7) + \frac{\zeta(5)\pi^2}{9} + \frac{17\zeta(3)\pi^4}{1260}$$

The next step is to find a closed expression for the alternate sums. We know from Apéry's 1979 article that it can be expressed in terms of $\zeta(3)$. If we omit the constant factor for short, we have:

$$\sum_{n=1}^{\infty} \frac{(-1)^{n+1}}{n^3 \binom{2n}{n}} = [\,\zeta(3)\,]$$

But unfortunately, there are no apparent patterns for $\zeta(5)$, $\zeta(7)$, etc. despite my numerous efforts for those sums and previous tests of D. H. Bailey and others are equally inconclusive so far. The mystery remains until someone finds the clue to all this.

An indication for further discoveries in that direction might be given with these remarks:

There are many representations of the same thing, for example,

$$\sum_{n=1}^{\infty} \frac{(-1)^{n+1}}{n \binom{2n}{n}} = [\,\psi\left(\frac{1}{5}\right), \psi\left(\frac{2}{5}\right), \psi\left(\frac{3}{5}\right), \psi\left(\frac{4}{5}\right)\,]$$

is in fact a linear combination of these numbers (when simplified).

$$\sum_{n=1}^{\infty} \frac{(-1)^{n+1}}{n \binom{2n}{n}} = [\ln 2\sqrt{5}, \ln(3+\sqrt{5})\sqrt{5}]$$

That is, a simple linear combination of $\ln(2)\sqrt{5}$ and $\ln(3+\sqrt{5})\sqrt{5}$. Which is also this number:

$$\frac{2}{5}\operatorname{arcsinh}\left(\frac{1}{2}\right)\sqrt{5}$$

This identity is well known and not difficult to find. So, the case $n^5$ should be simple and in terms of $\sqrt{5}$ and $\psi'''\left(\frac{1}{5}\right)$? This is strange since the other formulas were with $\pi\sqrt{3}$ only. It could be also that the $\sqrt{5}$ factor is an artefact. Normaly the sum for $n^5$ should be expressible with $\psi''''\left(\frac{k}{5}\right)$, $k = 1..4$ and/or $\psi''''\left(\frac{1}{3}\right)$ $\pi\sqrt{3}$ but numerical evidence shows that it is not the case.